\documentclass[12pt]{amsart}
\usepackage{amssymb}
\usepackage{verbatim}
\usepackage[usenames]{color}
\usepackage{hyperref}

\newtheorem{thm}{Theorem}
\newtheorem{prop}[thm]{Proposition}

\newtheorem{cor}[thm]{Corollary}

\theoremstyle{definition}
\newtheorem{df}[thm]{Definition}

\theoremstyle{remark}
\newtheorem{ex}[thm]{Example}
\newtheorem{rmk}[thm]{Remark}

\newenvironment{ls}{\begin{itemize}}{\end{itemize}}
\newenvironment{lsnum}{\begin{enumerate}}{\end{enumerate}}
\newenvironment{pf}{\begin{proof}}{\end{proof}}
\newcommand{\ger}[1]{\ensuremath{\mathfrak {#1}}}
\newcommand{\scr}[1]{\ensuremath{\mathcal {#1}}}

\newcommand{\bbb}[1]{\ensuremath{\mathbb {#1}}}

\renewcommand{\phi}{\varphi}

\newcommand{\notarrow}{\kern .42em\not\kern -.42em\longrightarrow}

\newcommand{\fe}{\ensuremath{\leq_{\text{fe}}}}

\newcommand{\hN}{{}^*{\bbb N}}
\newcommand{\hA}{{}^*\!A}
\newcommand{\hB}{{}^*\!B}

\newcommand{\hY}{{}^*Y}

\newcommand{\noprint}[1]{\relax}

\title{Finite Embeddability of Sets and Ultrafilters}
\author{Andreas Blass}
\address{Mathematics Department\\
University of Michigan\\
Ann Arbor, U.S.A.}
\email{ablass@umich.edu}
\thanks{Partially supported by NSF grant DMS-0653696}
\author{Mauro Di Nasso}
\address{Dipartimento di Matematica\\
Universit\`a di Pisa\\
Italy.}
\email{dinasso@dm.unipi.it}
\thanks{Partially supported by PRIN 2012 grant ``Logica, Modelli e Insiemi''.}

\begin{document}



\begin{abstract}
A set $A$ of natural numbers is finitely embeddable in another such
set $B$ if every finite subset of $A$ has a rightward translate that
is a subset of $B$.  This notion of finite embeddability arose in
combinatorial number theory, but in this paper we study it in its own
right.  We also study a related notion of finite embeddability of
ultrafilters on the natural numbers.  Among other results, we obtain
connections between finite embeddability and the algebraic and
topological structure of the Stone-\v Cech compactification of the
discrete space of natural numbers.  We also obtain connections with
nonstandard models of arithmetic.
\end{abstract}

\maketitle

\section{Introduction}  \label{intro}

The notion of finite embeddability of sets of natural numbers arose
naturally in combinatorial number theory \cite{diff}
(see also \cite{ruz} where the notion is implicitly used).
The present
paper is a study of the basic properties of this notion and a closely
related one in the realm of ultrafilters on the set \bbb N of natural
numbers. These notions have also been considered in
\cite{luperi-thesis}.  Additional information about them has appeared
in \cite{luperi-paper}, and more will appear
in a planned sequel to this paper.

\begin{df}[\cite{diff} \S 4]    \label{df-fe-set}
For $A,B\subseteq\bbb N$, we say that $A$ is \emph{finitely
embeddable} in $B$ and we write $A\fe B$ if each finite subset $F$
 of $A$ has a rightward translate $F+k$ included in $B$.
\end{df}

We use the standard notations $A+k=\{a+k:a\in A\}$ and $A-k=\{x\in\bbb
N:x+k\in A\}$ when $A\subseteq \bbb N$ and $k\in\bbb N$.  We also use
the standard conventions that the set \bbb N of natural numbers
contains 0 and that each natural number $n$ is identified with the set
$\{0,1,\dots,n-1\}$ of its predecessors.

\begin{rmk}
Our definition of finite embeddability differs from that in
\cite{diff} in that we work in \bbb N and use only rightward shifts.
The corresponding Definition~1.3 in \cite{diff} worked with subsets of
\bbb Z and allowed shifts in both directions.
\end{rmk}

\begin{df}      \label{df-fe-uf}
  For ultrafilters \scr U, \scr V on \bbb N, we say that \scr U is
  \emph{finitely embeddable} in \scr V and we write $\scr U\fe\scr V$ if, for
  each set $B\in\scr V$, there is some $A\in\scr U$ such that $A\fe
  B$.
\end{df}

It is clear that both of the relations \fe, one on subsets of \bbb N
and one on ultrafilters, are reflexive and transitive.

In Sections~\ref{sets} and \ref{ufs}, we study finite embeddability
primarily from a combinatorial point of view, with occasional mentions
of topological aspects.  The connection between finite embeddability
and nonstandard models, though crucial for the original motivation of
finite embeddability, has been postponed to Section~\ref{ns}, in
order to make most of our results accessible to readers unfamiliar
with nonstandard methods. On the other hand, readers who are
comfortable with nonstandard models can read Section~\ref{ns} without
first working through the preceding sections.

\section{Finite Embeddability of Sets of Natural
  Numbers}       \label{sets}

The following theorem summarizes some equivalent formulations
of finite embeddability of sets of natural numbers.
Additional equivalent characterizations in terms of nonstandard models will be
given in Section~\ref{ns}.

\begin{thm}            \label{fe-sets}
For any $A,B\subseteq\bbb N$, the following are equivalent.
\begin{lsnum}
\item $A\fe B$.
\item The family $\{B-a:a\in A\}$ has the finite intersection
  property.
\item There exists an ultrafilter \scr V on \bbb N such that $A$ is a
  subset of the ``leftward \scr V-shift'' of $B$, namely
\[
B-\scr V=\{x\in\bbb N:B-x\in\scr V\}.
\]
\item
 There exists an ultrafilter \scr V on \bbb N such that
 $A=B'-\scr V$ for some subset $B'$ of $B$.
\item The basic open sets $\overline A$ and $\overline B$ in the
  Stone-\v Cech compactification $\beta\bbb N$ satisfy $\overline
  A+\scr V\subseteq\overline B$ for some ultrafilter $\scr
  V\in\beta\bbb N$.
\item Some superset of $A$ is in the topological closure, in the power
  set $\scr P(\bbb N)$, of the set of leftward shifts $\{B-k:k\in\bbb
  N\}$ of $B$.
\item $A$ is in the topological closure of the set of leftward shifts
  of some subset $B'$ of $B$.
\end{lsnum}
\end{thm}

Before beginning the proof, we clarify the notation and terminology
used in this proposition, and we comment on some alternative ways to
view parts of it.

In item~(5), we identify the Stone-\v Cech compactification of \bbb N
with the set of ultrafilters on \bbb N, where natural numbers are
identified with the corresponding principal ultrafilters.  The basic
open sets are defined as $\overline A=\{\scr U\in\beta\bbb N:A\in\scr
U\}$, and this notation is justified by the fact that $\overline A$ is
also the closure in $\beta\bbb N$ of the set $A\subseteq\bbb
N\subseteq\beta\bbb N$.

The operation of addition of natural numbers is extended to $\beta\bbb
N$ by defining
\[
\scr U+\scr V=\{X\subseteq\bbb N:\{k:X-k\in\scr V\}\in\scr U\}.
\]
See \cite{hs} for extensive information about this operation (and its
analogs for other semigroups).  The notation $\overline A+\scr V$ in
item~(5) of the proposition means $\{\scr U+\scr V:\scr U\in\overline
A\}$.

In items~(6) and (7), the power set $\scr P(\bbb N)$ is to be
understood as topologized as the product $\{0,1\}^{\bbb N}$ via the
identification of subsets of \bbb N with their characteristic
functions; here $\{0,1\}$ is given the discrete topology.  So a small
neighborhood of a subset $A$ of \bbb N consists of those
$X\subseteq\bbb N$ that share a long initial segment with $A$.

The leftward shifts of a set, as used in items~(6) and (7), constitute
the orbit of that set under the transformation $X\mapsto X-1$ of the
space $\scr P(\bbb N)$.  So these two items can be reformulated in
terms of orbit closures for this transformation.

To the authors' knowledge,
the notion of ``leftward \scr V-shift'' in item~(3) was
first considered by Peter Krautzberger in his thesis \cite{krautz};
independently, Mathias Beiglb\"ock introduced it
for his ultrafilter proof of Jin's theorem \cite{mb}.

This item is a first indication
that finite embeddability is related to ultrafilters; that relation
plays an important role later in this paper.

\begin{pf}[Proof of Theorem~\ref{fe-sets}]

We shall first prove $(1)\iff(2)$, then $(3)\iff(5)$,
and finally the cycle
$(4)\implies(3)\implies(6)\implies(1)\implies(7)\implies(4)$.

  $(1)\iff(2)$ is immediate from the definitions, because $k$ is in a
  finite intersection $(B-a_1)\cap\dots(B-a_n)$ if and only if
  $\{a_1,\dots,a_n\}+k\subseteq B$.

$(3)\implies(5)$: If \scr V is as in (3), and if $\scr U\in\overline A$,
then $\{x\in\bbb N:B-x\in\scr V\}$, being a superset of $A$, is in
\scr U.  This means that $B\in\scr U+\scr V$, and so $\scr U+\scr
V\in\overline B$.  Since \scr U was arbitrary in $\overline A$, we
have (5).

$(5)\implies(3)$: Apply (5) to the principal ultrafilters $\scr U_a$
concentrated at points $a\in A$.  These are in $\overline A$, so we
have $B\in\scr U_a+\scr V$.  But this means that $B-a\in\scr V$, as
required for (3).

$(4)\implies(3)$: If \scr V and $B'$ are as in (4), then $A=B'-\scr
V\subseteq B-\scr V$.

$(3)\implies(6)$: Every leftward \scr V-shift $B-\scr V$, as defined
in (3), is the limit along \scr V of the leftward shifts $B-k$.  So the
set of all these leftward \scr V-shifts, for all \scr V, is the
topological closure mentioned in (6). Finally, the specific $B-\scr V$
asserted to exist in (3) serves as the superset of $A$ required in (6).

$(6)\implies(1)$: Let $F$ be any finite subset of $A$ and therefore
also of the superset $A'$ mentioned in (6).  By definition of the
topological closure, there must be a leftward shift $B-k$ whose
characteristic function agrees on $F$ with that of $A'$.  Therefore
$F\subseteq B-k$.

$(1)\implies(7)$: According to the assumption (1), for
each $n\in\bbb N$ there exists some $k\in\bbb N$ such that $(A\cap
n)+k\subseteq B$.  (Recall that we use the usual convention identifying a
natural number $n$ with its set of predecessors.)
If one and the same $k$ works for arbitrarily large $n$ (and therefore
for all $n$), then $A+k\subseteq B$.  Setting $B'=A+k$, we have (7),
because $B'\subseteq B$ and $A=B'-k$ is in the set of leftward shifts
of $B'$ (not just in its closure).

So assume from now on that, for each $k$, there is an upper bound on
the $n$'s for which $(A\cap n)+k\subseteq B$.  Then, for each
$n\in\bbb N$, there must exist arbitrarily large $k\in\bbb N$ with
$(A\cap n)+k\subseteq B$.  Indeed, let $n$ be given and consider any
$m\in\bbb N$; we shall find an appropriate (for $n$) $k\geq m$.  By
assumption, for each $k<m$, we have an upper bound $b_k$ on the $n$'s
for which this $k$ works; fix some $N$ larger than all of these $m$
bounds and larger than our given $n$.  There is a $k$ that works for
$N$ and thus also for our given $n$, and this $k$ must, by choice of
$N$, be $\geq m$.

Now define a subset $B'$ of $B$ as the union of finite pieces of the
form $(A\cap n)+k_n$, where the $k_n$'s are chosen (for $n\in\bbb
N$) inductively as follows.
Suppose we already have $k_j$ for all
$j<n$.  Then choose $k_n$ so that $k_n>k_j+j$ for all $j<n$ and
$(A\cap n)+k_n\subseteq B$.  This ensures that (the characteristic
function of) $A$ agrees on $[0,n)$ with that of a leftward shift
$B'-k_n$ of $B'$.  Since this happens for each $n$, $A$ is in the
closure of the set of leftward shifts of $B'$.

$(7)\implies(4)$: The closure of the set of leftward shifts $B'-k$
consists of exactly the sets $B'-\scr V$ for all ultrafilters \scr V
on \bbb N.
\end{pf}

\begin{rmk}
Peter Krautzberger independently proved the equivalence of items~(1),
(2), and (5) in Theorem~\ref{fe-sets}.    After learning of the second
author's proof that $\scr U\fe\scr U+\scr V$, he also proved
Theorem~\ref{closure-sums} below.
\end{rmk}

\begin{rmk}
  Suppose that, when we topologize $\scr P(\bbb N)\cong\{0,1\}^{\bbb
    N}$ as a product, we do not give $\{0,1\}$ the discrete topology
  but rather the Sierpi\'nski topology in which $\{1\}$ is open but
  $\{0\}$ is not.  Then $A\fe B$ if and only if $A$ is in the closure
  of the set of leftward shifts of $B$.  The reason is that, with this
  new topology, closed sets are automatically closed downward with
  respect to $\subseteq$.
\end{rmk}

As mentioned in the introduction, the notion of finite embeddability
was originally motivated by considerations from combinatorial number
theory.  To give an idea of this motivation, we list in the following
proposition some facts from \cite[Propositions~4.1 and 4.2]{diff};
for the relevant definitions, see also \cite{hs}.

\begin{prop}
Let $A$ and $B$ be sets of natural numbers.
  \begin{lsnum}
\item $A$ is maximal with respect to $\fe$ if and only if it is thick,
  i.e., it includes arbitrarily long intervals.
\item If $A\fe B$ and $A$ is piecewise syndetic, then $B$ is also
  piecewise syndetic.
\item If $A\fe B$ and $A$ contains a $k$-term arithmetic progression,
  then also $B$ contains a $k$-term arithmetic progression.
\item If $A\fe B$ then the upper Banach densities satisfy $BD(A)\leq
  BD(B)$.
\item If $A\fe B$ then $A-A\subseteq B-B$.
\item If $A\fe B$ then $\bigcap_{t\in G}(A-t)\fe\bigcap_{t\in G}(B-t)$
  for every finite $G\subseteq\bbb N$.
  \end{lsnum}
\end{prop}

We remark in connection with item~(2) that ``piecewise'' is
essential there.  The property of being syndetic is not in general preserved
upward under $\fe$.  Intuitively, this is because the definition of
$A\fe B$ does nothing to prevent the occurrence of long gaps in $B$.
For similar reasons, ``upper'' is essential in item~(4); asymptotic
density is not monotone with respect to $\fe$.

\section{Finite Embeddability of Ultrafilters}  \label{ufs}

The results of the preceding section exhibit some connections between
finite embeddability and the additive structure of ultrafilters.  In
the present section, we first reformulate those connections in some
corollaries and a theorem, and then we establish some additional connections.

\begin{df}
  Let \scr U be an ultrafilter on \bbb N.  A set $B\subseteq\bbb N$ is
  \scr U-\emph{rich} if there is an $A\in\scr U$ with $A\fe B$.
\end{df}

Notice that Definition~\ref{df-fe-uf} of finite embeddability of
ultrafilters says that $\scr U\fe\scr V$ if and only if every set in
\scr V is \scr U-rich.

\begin{cor}
Let \scr U be an ultrafilter on \bbb N.  A set $B$ is \scr U-rich if
and only if $B\in\scr U+\scr W$ for some ultrafilter \scr W on \bbb N.
\end{cor}

\begin{pf}
The equivalence of (1) and (3) in Proposition~\ref{fe-sets} shows that
$B$ is \scr U-rich if and only if there is an ultrafilter \scr W on
\bbb N such that $B-\scr W\in\scr U$.  But this is precisely what is
required for $B\in\scr U+\scr W$.
\end{pf}

\begin{cor}     \label{part-reg}
  For any ultrafilter \scr U, the family of \scr U-rich sets is
  parti\-tion-regular.  That is, if a \scr U-rich set is partitioned
  into finitely many pieces, then at least one of the pieces is
  \scr U-rich.
\end{cor}

\begin{pf}
  The preceding corollary shows that the family of \scr U-rich sets is
  the union of a collection of ultrafilters, and it is easy to see
  that any such union is partition-regular.
\end{pf}

\begin{thm}     \label{closure-sums}
Let \scr U and \scr V be ultrafilters on \bbb N.  Then $\scr U\fe\scr
V$ if and only if \scr V is in the closure, in $\beta\bbb N$, of the
set of sums $\{\scr U+\scr W:\scr W\in\beta\bbb N\}$.
\end{thm}

\begin{pf}
Each of the following statements is clearly equivalent to the next.
\begin{ls}
\item $\scr U\fe\scr V$.
\item Every $B\in\scr V$ is \scr U-rich.
\item Every $B\in\scr V$ is in $\scr U+\scr W$ for some ultrafilter
  \scr W on \bbb N.
\item Every basic neighborhood $\overline B$ of \scr V in $\beta\bbb
  N$ contains a sum $\scr U+\scr W$.
\item \scr V is in the closure of the set $\{\scr U+\scr W:\scr
  W\in\beta\bbb N\}$ of sums. \qedhere
\end{ls}
\end{pf}

The preceding theorem
implies, in particular, that the upward cone
$\{\scr V:\scr U\fe\scr V\}$ determined by any ultrafilter \scr U in
the $\fe$ ordering is a closed subset of $\beta\bbb N$.

This fact can also be seen by inspecting the form of the definition of
$\fe$.  To say that \scr V is in this cone is to say ``For all
$B\in\scr V$, $\Phi(B)$'' where the statement $\Phi(B)$ doesn't
mention \scr V.  Any set of ultrafilters with a definition of this
form is closed.  Furthermore, it can be expressed as the set of all
ultrafilters extending the filter generated by the complements of all
sets $B$ that don't satisfy $\Phi(B)$.

In the case at hand, we have that the upward cone $\{\scr V:\scr
U\fe\scr V\}$ consists of all extensions of the filter $\scr F_{\scr
  U}$ generated by those sets whose complements are not \scr U-rich.
In view of Corollary~\ref{part-reg}, this description can be
simplified, because the sets whose complements are not \scr U-rich
already constitute a filter $\scr F_{\scr U}$; there is no need to
form the filter they generate.

The theorem provides another description of $\scr F_{\scr U}$ as an
intersection of ultrafilter sums,
\[
\scr F_{\scr U}=\bigcap_{\scr W\in\beta\bbb N}(\scr U+\scr W).
\]

The preceding results have related finite embeddability of \scr U with
sums in which \scr U appears as the left summand.  Since addition is
not commutative in $\beta\bbb N$, the question arises whether there
are similar results with \scr U as the right summand.  Most of the
preceding results do not have such analogs, but the basic fact (a
special case of Theorem~\ref{closure-sums}) that $\scr U\fe\scr
U+\scr W$ does, and in fact we get a somewhat stronger conclusion.

\begin{thm}    \label{left-sum}
For any ultrafilters \scr V and \scr W on \bbb N,
the following are equivalent:
\begin{lsnum}
\item
Each set in \scr W includes a rightward translate of a whole set
from \scr V;
\item
$\scr W=\scr U+\scr V$ for a suitable \scr U.
\end{lsnum}
\end{thm}

\begin{pf}
$(1)\implies(2)$: The family $\scr F=\{B-\scr V:B\in\scr W\}$
of the leftward \scr V-shifts of sets in \scr W
has the finite intersection property. Indeed, for each $B\in\scr W$,
by hypothesis one has that $A+k\subseteq B$ for suitable $A\in\scr V$ and
$k\in\bbb N$, and so $k\in B-\scr V\ne\emptyset$.
Moreover, given finitely many $B_1,\ldots,B_n\in\scr W$,
the intersection $\bigcap_{i=1}^n(B_i-\scr V)=\left(\bigcap_{i=1}^n B_i\right)-\scr V\in\scr F$,
and therefore it is nonempty.
So we can extend \scr F to an ultrafilter \scr U.
For every
$B\in\scr W$ one has that $B-\scr V\in\scr F\subseteq\scr U$,
so $B\in\scr U+\scr V$, and the equality $\scr W=\scr U+\scr V$ follows.

$(2)\implies(1)$: The definition of $\scr U+\scr V$ says that if
$B\in\scr U+\scr V$ then for some $k$, in fact for \scr U-almost all
$k$, we have $B-k\in\scr V$.  So $B-k$ is a set in \scr V with a
rightward translate, by $k$, included in $B$.
\end{pf}

Remark that $(1)$ in Theorem~\ref{left-sum} above is a strictly
stronger property than $\scr V\fe\scr W$.

\begin{cor}
The ordering $\fe$ of ultrafilters on \bbb N is upward directed.
\end{cor}

\begin{pf}
Any \scr U and \scr V have $\scr U+\scr V$ as an upper bound.
\end{pf}

In view of the characterization of the upward cone of \scr U as the
(topological) closure of $\{\scr U+\scr W:\scr W\in\beta\bbb N\}$, it
is tempting to consider the closure of $\{\scr W+\scr U:\scr
W\in\beta\bbb N\}$ also, but this leads to nothing interesting.
Indeed, $\{\scr W+\scr U:\scr W\in\beta\bbb N\}$ is already closed,
being the image of the compact set $\beta\bbb N$ under the continuous
function $\scr W\to\scr W+\scr U$.  (Recall that addition in
$\beta\bbb N$ is a continuous function of the left summand when the
right summand is fixed, but not vice versa.)

Let us summarize the preceding results, in the language of ideals of
the semigroup $(\beta\bbb N,+)$.  Recall that $\{\scr U+\scr W:\scr
W\in\beta\bbb N\}$ and $\{\scr W+\scr U:\scr W\in\beta\bbb N\}$ are,
respectively, the right and left ideals generated by \scr U.

\begin{cor}
For any $\scr U\in\beta\bbb N$, the upward cone $\{\scr V:\scr
U\fe\scr V\}$ is a closed, two-sided ideal in $\beta\bbb N$.  It is
the smallest closed right ideal containing \scr U, and therefore it is
also the smallest closed two-sided ideal containing \scr U.
\end{cor}

The fact that the closure of a right ideal is a two-sided ideal is
a general property of Stone-\v Cech compactifications of discrete
commutative semigroups.  In fact, it holds in even greater generality;
see Theorem~2.19(a) of \cite{hs}.

For completeness, we mention that we cannot replace ``right'' by
``left'' in this corollary.  The smallest left ideal containing \scr U,
namely $\beta\bbb N+\scr U$, is already closed and, as the following
example shows, it need not contain the right ideal generated by \scr
U.

\begin{ex}
  Recall that the sets $X\subseteq\bbb N$ whose density (i.e., the
  limit of $|X\cap n|/n$ as $n\to\infty$) exists and is zero
  constitute a proper ideal, so there are ultrafilters \scr U that
  contain no such $X$.  The set of all such ultrafilters \scr U is
  easily seen to be a left-ideal $\Delta$ (see
  \cite[Theorem~6.79]{hs}).  So, if we fix one such ultrafilter \scr
  U, then the left ideal it generates is included in $\Delta$.

 We claim that, in contrast, the right ideal generated by \scr U
  contains an ultrafilter $\scr U+\scr V$ that contains a set of
  density zero.  To see this, let \scr V be any non-principal
  ultrafilter on \bbb N that contains the set $P=\{2^m:m\in\bbb N\}$,
  and consider the set
\[
Q=\{2^m+k:m,k\in\bbb N\text{ and }k<m\}.
\]
Observe that, for every fixed $k$, the set $Q-k$ contains all but
finitely many elements of $P$ and is therefore in \scr V.  Thus,
$Q\in\scr U+\scr V$.  (Indeed, $Q\in\scr U'+\scr V$ for all
ultrafilters $\scr U'$ on \bbb N.) Finally, to see that $Q$ has
density 0, notice that the function $n\mapsto |A\cap n|/n$ achieves
its local maxima at the numbers of the form $n=2^m+m$, where its
values are $(m^2+m)/(2^{m+1}+2m)$.  These values tend to zero as
$m\to\infty$, so $Q$ has density zero.  (The essential properties of
$Q$ used here are merely that it has density zero but is thick.
Recall that thickness means that every finite set has a rightward
translate included in $Q$, or equivalently that there is an
ultrafilter \scr V that contains all the sets $Q-k$.)
\end{ex}

\section{Nonstandard Models}    \label{ns}

In this section we assume the reader to be familiar with the basics
of nonstandard analysis, and in particular with the notions of
\emph{hyper-extension} $\hA$ (or nonstandard extension) of a set $A$,
and with the \emph{transfer principle}. Most notably,
we will focus on the set of \emph{hypernatural numbers} $\hN$,
and we will consider the nonstandard characterization of topology.
We will work in a nonstandard framework where
the $\ger c^+$-\emph{enlarging property} holds,
namely the property that if a family $\ger X$ of cardinality
at most the continuum has the finite intersection property,
then the intersection $\bigcap_{X\in\ger X}{}^*X$
is nonempty. (This is a weaker property than
$\ger c^+$-saturation.)
We refer to \S 4.4 of \cite{ck} for the foundations,
and to \cite{ma} for all nonstandard notions
and results used in this section.

Finite embeddability has a suggestive nonstandard characterization,
which was in fact the original motivation for this research.
Precisely, $A\fe B$ means that a (possibly infinite) rightward translation
of $A$ is included in the hyper-extension of $B$.

\begin{prop}[\S 4 of \cite{diff}] \label{fe-nonstandard}
$A\fe B$ if and only if $\mu+A\subseteq\hB$ for some $\mu\in\hN$.
\end{prop}

\begin{pf}
By the definition of $A\fe B$, it directly follows that
the family $\{B-a:a\in A\}$ satisfies the finite intersection
property. Then, by the enlarging property,
we can pick $\mu\in\bigcap_{a\in A}{}^*(B-a)$.
Since ${}^*(B-a)=\hB-a$, this means that $\mu+A\subseteq\hB$.

Conversely, given a finite subset $\{a_1,\ldots,a_k\}\subseteq A$,
the element $\mu$ witnesses that the following property holds:
$$\exists x\in\hN\ x+a_1\in\hB\land\ldots\land x+a_k\in\hB.$$
Then, by transfer, we obtain the existence of an element $x\in\bbb N$
such that $x+a_i\in B$ for all $i=1,\ldots,k$.
\end{pf}

There is a canonical way of associating an ultrafilter on $\bbb N$ to
each hypernatural number.

\begin{df}
The \emph{ultrafilter generated} by $\alpha\in\hN$ is the family
$$\scr U_\alpha\ =\ \{A\subseteq\bbb N:\alpha\in\hA\}$$
\end{df}

It is readily checked that $\scr U_\alpha$ is actually an ultrafilter,
and that $\scr U_\alpha$ is non-principal if and only if
$\alpha\notin\bbb N$, i.e., $\alpha$ is infinite.
Notice that every ultrafilter on $\bbb N$ is a family with
the cardinality of the continuum that satisfies
the finite intersection property.
So, by the $\ger c^+$-enlarging property, we obtain
that each ultrafilter is generated by some
(actually, by $\ger c^+$ many)
hypernatural numbers; in consequence:
$$\beta\bbb N\ =\ \{\scr U_\mu:\mu\in\hN\}.$$

For $B\subseteq\bbb N$ and $\mu\in\hN$, we consider
the set of elements in $\hB$ that are
placed at finite distance from $\mu$ on the right side:
$$B_\mu\ =\ \{n\in\bbb N: \mu+n\in\hB\}\ =\ (\hB-\mu)\cap\bbb N.$$
Trivially $B_\mu\fe B$ because $\mu+B_\mu\subseteq\hB$.
Notice that $B_\mu$ is the leftward $\scr U_\mu$-shift of $B$; indeed
$n\in B-\scr U_\mu\Leftrightarrow B-n\in\scr U_\mu\Leftrightarrow
\mu\in{}^*(B-n)\Leftrightarrow \mu+n\in\hB$. In consequence:
$$A\in\scr U_\alpha+\scr U_\beta\iff
A_\beta=A-\scr U_\beta\in\scr U_\alpha\iff
\alpha\in\hA_\beta.$$

A nice nonstandard characterization can also be given
of the topology on $\scr P(\bbb N)$
considered in Section \ref{sets}.

\begin{prop} \label{topology-nonstandard}
$A$ is in the topological closure, in the power
set $\scr P(\bbb N)$, of the set of leftward shifts
$\{B-k:k\in\bbb  N\}$ of $B$ if and only if
$A=B_\mu$ for some $\mu\in\hN$.
\end{prop}

\begin{pf}
Recall that if $X$ is a topological space of character $\kappa$,
then in any nonstandard model with the $\kappa^+$-enlarging property,
the following characterization holds:
``An element $x$ belongs to the closure of a set
$Y\subseteq X$ if and only if $m(x)\cap\hY\ne\emptyset$,'' where
$$m(x)\ =\ \bigcap\{{}^*I:I \text{ neighborhood of }x\}$$
is the monad of $x$. (See \emph{e.g.} \S 3.1 of \cite{ma}.)
In our topological space, a base of neighborhoods
of a set $A\subseteq\bbb N$ is obtained by taking
all finite intersections of sets
of the form $I_{A,n}=\{B\in\scr P(\bbb N):n\in A\Leftrightarrow n\in B\}$.
So, we have that:
$$m(A)\ =\ \bigcap_{n\in\bbb N}{}^*I_{A,n}\ =\
\{X\in{}^*{\scr P}(\bbb N):X\cap\bbb N=A\}.$$
By the above nonstandard characterization,
$A$ is in the closure of the set $\{B-k:k\in\bbb  N\}$
if and only if there exists $\mu\in\hN$ such that
$\hB-\mu\in m(A)$, \emph{i.e.} $B_\mu=(\hB-\mu)\cap\bbb N=A$.
\end{pf}

The above notions and characterizations can be used to give
nonstandard proofs of all the results presented in this paper.
Below, we consider in detail the main results (named here Theorems A and B)
and leave the other proofs as exercises for the interested reader.

\medskip
\noindent
\textbf{Theorem A.} [Theorem \ref{fe-sets}]
\\
\emph{
For any $A,B\subseteq\bbb N$, the following are equivalent.
\begin{lsnum}
\item $A\fe B$.
\item The family $\{B-a:a\in A\}$ has the finite intersection
  property.
\item There exists an ultrafilter \scr V on \bbb N such that $A$ is a
  subset of the ``leftward \scr V-shift'' of $B$, namely
\[
B-\scr V=\{x\in\bbb N:B-x\in\scr V\}.
\]
\item
 There exists an ultrafilter \scr V on \bbb N such that
 $A=B'-\scr V$ for some subset $B'$ of $B$.
\item The basic open sets $\overline A$ and $\overline B$ in the
  Stone-\v Cech compactification $\beta\bbb N$ satisfy $\overline
  A+\scr V\subseteq\overline B$ for some ultrafilter $\scr
  V\in\beta\bbb N$.
\item Some superset of $A$ is in the topological closure, in the power
  set $\scr P(\bbb N)$, of the set of leftward shifts $\{B-k:k\in\bbb
  N\}$ of $B$.
\item $A$ is in the topological closure of the set of leftward shifts
  of some subset $B'$ of $B$.
\item
There exists $\mu\in\hN$ such that $A\subseteq B_\mu$.
\item
$A'=B_\mu$ for some superset $A'$ of $A$ and some $\mu\in\hN$.
\item
$A=B'_\mu$ for some subset $B'$ of $B$ and some $\mu\in\hN$.
\end{lsnum}}

\begin{pf}[Nonstandard proof]
We first reduce to the ``nonstandard" items $(8)$,
$(9)$ and $(10)$.

$(1)\iff (8)$. It is Proposition \ref{fe-nonstandard}.

$(2)\iff (8)$. Recall that a family has the finite intersection property
if and only if it is included in an ultrafilter. So, item (2) is equivalent
to the existence of a point $\mu\in\hN$ such that
$\{B-a:a\in A\}\subseteq\scr U_\mu$. Now,
$B-a\in\scr U_\mu\Leftrightarrow \mu\in{}^*(B-a)=\hB-a\Leftrightarrow a\in B_\mu$,
and therefore $\{B-a:a\in A\}\subseteq\scr U_\mu$ is equivalent to $A\subseteq B_\mu$.

$(3)\iff (8)$. Item (3) says that $A\subseteq B-\scr U_\mu$ for
some $\mu\in\hN$. Then recall that $B-\scr U_\mu=B_\mu$.

$(4)\iff (10)$. Condition (4) says that $A=B'-\scr U_\mu=B'_\mu$ for some
$B'\subseteq B$ and some $\mu\in\hN$.

$(5)\iff (8)$. Notice that
$\overline{A}=\{\scr U_\alpha:\alpha\in\hA\}$ and
$\overline{B}=\{\scr U_\beta:\beta\in\hB\}$.
So, item (5) says that there exists $\mu\in\hN$ such that
$B\in\scr U_\alpha+\scr U_\mu$ for all $\alpha\in\hA$.
Now recall that $B\in\scr U_\alpha+\scr U_\mu$
if and only if $\alpha\in\hB_\mu$.
Then the above property is equivalent to $\hA\subseteq\hB_\mu$,
which in turn is equivalent to $A\subseteq B_\mu$, by transfer.

$(6)\iff (9)$ and $(7)\iff (10)$. They both directly follow
from Proposition \ref{topology-nonstandard}.

We are left to prove the equivalence of the three ``nonstandard" conditions.

$(8)\iff (9)$. Trivial.

$(10)\implies (8)$. By transfer, $B\subseteq B'\Leftrightarrow \hB\subseteq\hB'$,
and hence $B\subseteq B'\Rightarrow B_\mu\subseteq B'_\mu$ for all $\mu\in\hN$.

$(8)\implies (10)$. Assume first that $\mu=k\in\bbb N$ is finite.
By the hypothesis, $A\subseteq B_k=B-k$.
If we let $B'=A+k$ then $B'\subseteq B$ and $A=B'_k$.

Now let $\mu\in\hN$ be infinite.
Denote by $A_n=\{a_1<\ldots<a_n\}$ the set of the first $n$ elements of $A$.
Notice that every set $\Lambda_n=\{x\in\bbb N:x+A_n\subseteq B\}$ is infinite,
since the infinite number $\mu\in{}^*\Lambda_n$.
We now inductively define a sequence of numbers $\{x_n\}$
and a sequence $\{B_n\}$ of finite subsets of $B$ as follows.
\begin{ls}
\item Let $x_1\in\Lambda_1$ and let $B_1=\{x_1+a_1\}=x_1+A_1$.
\item At the inductive step $n>1$, pick $x_n\in\Lambda_n$ such
that $x_n>x_{n-1}+a_{n-1}$ and let $B_n=x_n+A_n\subseteq [x_n,x_n+a_n]$.
\item Define $B'=\bigcup_{n\in\bbb N}B_n\subseteq B$.
\end{ls}
The sets $B_n$ are pairwise disjoint and so, for every $n$,
$$x_n+A_n\ =\ x_n+(A\cap[0,a_n])\ =\ B'\cap[x_n,x_n+a_n],$$
\emph{i.e.} $A\cap[0,a_n]=(B'-x_n)\cap [0,a_n]$.
Then, by transfer, for all $N\in\hN$,
$$\hA\cap[0,a_N]\ \ =(\hB'-x_N)\cap[0,a_N].$$
If we pick $\nu=x_N$ for an infinite $N$, then all finite numbers
$\bbb N\subseteq[0,a_N]$, and hence $A=\hA\cap\bbb N=(\hB'-\nu)\cap\bbb N=B'_\nu$.
\end{pf}

\medskip
\noindent
\textbf{Theorem B.} [Theorem \ref{closure-sums}]
\\
\emph{Let \scr U and \scr V be ultrafilters on \bbb N.  Then $\scr U\fe\scr
V$ if and only if \scr V is in the closure, in $\beta\bbb N$, of the
set of sums $\{\scr U+\scr W:\scr W\in\beta\bbb N\}$.}

\begin{pf}[Nonstandard proof]
Pick $\alpha,\beta\in\hN$ such that
$\scr U=\scr U_\alpha$ and $\scr V=\scr U_\beta$.
Assume first that $\scr U_\alpha\fe\scr U_\beta$.
We want to show that
for every $B\in\scr U_\beta$ there exists $\mu\in\hN$
such that $B\in\scr U_\alpha+\scr U_\mu$.
By the hypothesis we can pick
$A\in\scr U_\alpha$ with $A\fe B$.
By the nonstandard characterization, this means that
there exists $\mu\in\hN$ such that $A\subseteq B_\mu$.
But then $B_\mu\in\scr U_\alpha$, \emph{i.e.} $\alpha\in\hB_\mu$,
and we conclude that $B\in\scr U_\alpha+\scr U_\mu$.

Conversely, if $\scr U_\beta$ is in the closure
of $\{\scr U_\alpha+\scr U_\mu:\mu\in\hN\}$,
then for every $B\in\scr U_\beta$
there exists $\mu$ such that $B\in\scr U_\alpha+\scr U_\mu$,
\emph{i.e.} $\alpha\in\hB_\mu$.
But then we have found a set $B_\mu\in\scr U_\alpha$
with $B_\mu\fe B$, as desired.
\end{pf}

\section{Questions}     \label{questions}

\begin{ls}
\item
Under which conditions does the following implication hold?
\\
$A\fe B\ \text{and}\ A'\fe B'\implies (A-A')\fe(B-B')$.
\item
For $A\subseteq\bbb N$, is there a neat combinatorial description
of the equivalence classes $[A]=\{B:A\fe B\land B\fe A\}$?
And, for $\scr U\in\beta\bbb N$, of the equivalence classes
$[\scr U]=\{\scr V:\scr U\fe\scr V\land \scr V\fe \scr U\}$?
\item
If we modify the definition of $\fe$ so that
trivial right translations by $0$ are not permitted,
can we find a neat combinatorial description of
the sets $A$ such that $A\fe A$?
And of the ultrafilters \scr U such that $\scr U\fe\scr U$?
\end{ls}

\begin{rmk}
  Proposition~4.2(6) of \cite{diff} looks as if it answers the first of
  these questions, but its proof relies on the definition in
  \cite{diff} of finite embeddability, which allows both leftward and
  rightward shifts.  Isaac Goldbring has pointed out, in a private
  communication, that the same proof works in our setting, using only
  rightward shifts, under the additional hypothesis that, for each
  finite $F\subseteq A$, there are infinitely many $k\in\bbb N$ such that
  $F+k\subseteq B$.  Goldbring suggests calling this hypothesis
  \emph{proper finite embeddability} of $A$ in $B$, and notes that it
  is an intermediate property between our finite embeddability and
  dense embeddability as defined in \cite{diff}.  Note that, as shown
  in the proof of $(1)\implies(7)$ in Theorem~\ref{fe-sets} above, the
  only way $A$ can be finitely embeddable but not properly finitely
  embeddable in $B$ is that there is some $k\in\bbb N$ with
  $A+k\subseteq B$.

 It remains open
  whether proper finite embeddability is the exact condition for the
  first question or whether some weaker condition might suffice.
\end{rmk}

\end{document}